\documentclass[10pt]{article}
\usepackage{amssymb, amsmath, amsthm, color,graphicx, amsfonts,a4}
\usepackage{amscd}

\definecolor{Red}{cmyk}{0,1,1,0}

\definecolor{verde}{cmyk}{1,0,1,0}

\definecolor{azul}{cmyk}{1,1,0,0}

\numberwithin{equation}{section}

\newcommand{\E}{\mathbb{E}}
\newcommand{\N}{\mathbb{N}}

\newcommand{\Z}{\mathbb{Z}}
\renewcommand{\P}{\mathbb{P}}

\newcommand{\V}{\mathbb{V}}

\newcommand{\be}{\begin{equation}}
\newcommand{\ee}{\end{equation}}

\newtheorem{teorema}{Theorem}

\newtheorem{lema}{Lemma}

\newtheorem{conjectura}{Conjecture}

\begin{document}
\title{Critical Point and Percolation Probability in a Long Range Site Percolation Model on $\Z^d$}
\author{Bernardo N.B. de Lima\footnote{ Departamento de Matem{\'a}tica, Universidade Federal de Minas Gerais, Av. Ant\^onio
Carlos 6627 C.P. 702 CEP30123-970 Belo Horizonte-MG, Brazil} ,
R\'emy Sanchis$^{*}$, Roger W.C. Silva\footnote{ Departamento de
Estat\'\i stica, Universidade Federal de Minas Gerais, Av.
Ant\^onio Carlos 6627 C.P. 702 CEP30123-970 Belo Horizonte-MG,
Brazil and Departamento de
Matem\'atica, Universidade Federal de Ouro Preto, Rua Diogo de Vasconcelos 122 CEP35400-000 Ouro Preto-MG}}
\maketitle
\begin{abstract}
Consider an independent site percolation model with parameter $p \in
(0,1)$ on $\Z^d,\ d\geq 2$ where there are only nearest neighbor
bonds and long range bonds of length $k$ parallel to each coordinate
axis. We show that the percolation threshold of such model converges
to $p_c(\Z^{2d})$ when $k$ goes to infinity, the percolation
threshold for ordinary (nearest neighbour) percolation on $\Z^{2d}$.
We also generalize this result for models whose long range bonds
have several lengths.
\end{abstract}
{\footnotesize Keywords: long range percolation; percolation threshold }

\section{Introduction and Notation}

Let $G=(\V,\E)$ be a graph with a countably infinite vertex set
$\V$. Consider the Bernoulli site percolation model on $G$; to each
site $v\in\V$ we associate a Bernoulli random variable $X(v)$, which
takes the values 1 and 0 with probability $p$ and $1-p$
respectively. This can be done considering the probability space
$(\Omega, \mathcal{F},\P_p)$, where $\Omega=\{0,1\}^{\V}$,
$\mathcal{F}$ is the $\sigma$-algebra generated by the cylinder sets
in $\Omega$ and $\P_p=\prod_{v\in\V}\mu(v)$ is the product of
Bernoulli measures with parameter $p$, in which the configurations
$\{X(v),v \in \mathbb{V}\}$ take place. We denote a typical element
of $\Omega$ by $\omega$. When $X(v)=1$ (respectively, $X(v)=0$) we
say that $v$ is ``open'' (respectively, ``closed'').

Given two vertices $v$ and $u$, we say that $v$ and $u$ are
connected in the configuration $\omega$ if there exists a finite
path $\langle v=v_0,v_1,\dots,v_n=u\rangle$ of open vertices in
$\V$, such that $v_i\neq v_j,\ \forall i\neq j$ and $\langle
v_i,v_{i+1} \rangle$ belongs to $\E$ for all $i=0,1,\dots,n-1$. We
will use the short notation $\{v\leftrightarrow u\}$ to denote the
set of configurations where $u$ and $v$ are connected.

Given the vertex $v$, the cluster of $v$ in the configuration
$\omega$ is the set $C_v(\omega)=\{u\in \V; v\leftrightarrow u \mbox{ on } \omega\}$. We say that the vertex $v$
percolates when the cardinality of $C_v(\omega)$ is infinite; we
will use the following standard notation
$\{v\leftrightarrow\infty\}\equiv\{\omega\in \Omega;
\#C_v(\omega)=\infty\}$. Fixed some vertex $v$, we define the
percolation probability of the vertex $v$ as the function
$\theta_v(p):[0,1]\mapsto [0,1]$ with
$\theta_v(p)=\P_p(v\leftrightarrow\infty)$.

From now on, the vertex set $\V$ will be $\Z^d,\ d\geq 2$ and for
each positive integer $k$ define $$\E_k=\{\langle
(v_1,\dots,v_d),(u_1,\dots, u_d)\rangle\in \V\times\V; \exists !
i\in\{1,\dots,d\}\mbox{ such that }|v_i-u_i|=k\mbox{ and }v_j=
u_j,\forall j\neq i\}.$$ Let's define the graph $G^k=(\V,
\E_1\cup\E_k)$ that is, $G^k$ is $\Z^d$ equipped with nearest
neighbor bonds and long range bonds with length $k$ parallel to each
coordinate axis. Observe that $G^k$ is a transitive graph, hence the
function $\theta_v(p)$ does not depend on $v$ and we write only
$\theta^k(p)$ to denote $\P_p(0\leftrightarrow\infty)$ for any transitive graph.

The simplest version of our main result (see Theorem \ref{prin}
below) states that $p_c(G^k)$ converges to $p_{c}(\mathbb{Z}^{2d})$,
when $k$ goes to infinity. The main motivation to study this
question is that we believe that the Conjecture \ref{conjec} stated
below can shed some light on the truncation problem for long range
percolation. This problem, proposed by E. Andjel, is the following:

On $\Z^d,\ d\geq 2$, consider the complete graph $G=(\Z^d,\E)$, that
is for all $v,u\in\Z^d$ we have that $\langle u,v \rangle\in\E$. For
each bond $\langle u,v \rangle\in\E$ we define its length as
$\|u-v\|_1$. Given a sequence $(p_n\in[0,1), n\in\N)$ consider an
independent bond percolation model where each bond whose length is
$n$ will be open with probability $p_n$. Assume that
$\sum_{n\in\N}p_n=\infty$, by Borel-Cantelli's Lemma the origin will
percolate to infinity with probability 1. The general and still open
truncation question is the following: is it true that there exists
some sufficiently large but finite integer K, such that the origin
in the truncated processes, obtained deleting (or closing) all long
range bonds whose length are bigger than $K$, still percolates to
infinity with positive probability?


The general question is still open, see \cite{FL} for a more
detailed discussion.


\section{The main result}

Given a positive integer $n$, define the $n$-vector
$\overrightarrow{k}=(k_1,\dots,k_n)$, where $k_i\in\{2,3,\dots
\},\forall i=1,\dots,n$. We define the graph
$G^{\overrightarrow{k}}$ as $(\Z^d,\E_1\cup(\cup_{i=1}^n
\E_{k_1\times\dots\times k_i}))$. Observe that when $n=1$ and
$k_1=k$, the graph $G^{\overrightarrow{k}}$ is the graph $G^k$
defined above. That is, $G^{\overrightarrow{k}}$ is $\Z^d$ decorated
with all bonds parallel to each coordinate axis with lengths
$1,k_1,k_1\times k_2,\dots, k_1\times k_2\times\dots\times k_n.$

From now on, we will use the notation $S^{\overrightarrow{k}}$ to
denote the $d(n+1)-$dimensional slab graph where the vertex set is
$(\Z\times\prod_{i=1}^n \{0,1,\dots,k_{n-i+1}-1\})^d$ and
$S^{\overrightarrow{k}}$ is equipped with only nearest
neighbor bonds.

The aim of this note is to prove that the percolation function of
the graph $G^{\overrightarrow{k}}$ is bounded between the
percolation functions of the slab $S^{\overrightarrow{k}}$ and of
$\mathbb{Z}^{d(n+1)}$. More precisely, we have the following Lemmas:

\begin{lema}\label{inferior} For any $p\in [0,1],\ \theta_v^{S^{\overrightarrow{k}}}(p)\leq\theta^{G^{\overrightarrow{k}}}(p), \forall v\in\V(S^{\overrightarrow{k}}).$
\end{lema}
\begin{lema}\label{superior} For any $p\in [0,1],\ \theta^{G^{\overrightarrow{k}}}(p)\leq\theta^{\Z^{d(n+1)}}(p).$
\end{lema}

The proof of these Lemmas will be given in the next section. Combining Lemmas \ref{inferior} and \ref{superior} it is straightforward that
$$\theta_v^{S^{\overrightarrow{k}}}(p)\leq\theta^{G^{\overrightarrow{k}}}(p)\leq\theta^{\Z^{d(n+1)}}(p),\forall v\in\V(S^{\overrightarrow{k}}).$$
Moreover,

\begin{equation}\label{sanduba}p_c(\Z^{d(n+1)})\leq
p_c(G^{\overrightarrow{k}})\leq p_c(S^{\overrightarrow{k}}).
\end{equation}

Using Theorem A of \cite{GM}, we have that
\begin{equation}\label{slab}\lim_{k_i\rightarrow\infty,\forall
i}p_c(S^{\overrightarrow{k}})=p_c(\Z^{d(n+1)}).
\end{equation}
Then combining (\ref{sanduba}) and (\ref{slab}), we can conclude that
$$\lim_{k_i\rightarrow\infty,\forall
i}p_c(G^{\overrightarrow{k}})=p_c(\Z^{d(n+1)}).$$

We have just proved that main result of this paper.

\begin{teorema}\label{prin}
Let $p_c(\Z^{d(n+1)})$ be the ordinary site percolation threshold
for $\Z^{d(n+1)}$ with nearest neighbor connections, then
\begin{equation}\label{eqprin}\lim_{k_i\rightarrow\infty,\forall
i}p_c(G^{\vec{k}})=p_c(\Z^{d(n+1)}).\end{equation}
\end{teorema}

\section{Proofs of the Lemmas}

{\it Proof of Lemma \ref{inferior}}

From the graph $G^{\overrightarrow{k}}$ we define the graph
$F^{\vec{k}}$ deleting some bonds in $G^{\overrightarrow{k}}$. More
precisely, $F^{\overrightarrow{k}}$ is the graph
$(\Z^d,(\E_1\cup(\cup_{i=1}^n \E_{k_1\times\dots\times
k_i}))-\cup_{i=1}^nR_i)$, where $R_i=R_i(\overrightarrow k)$ is the
set of deleted bonds with length $k_1\times\dots\times k_{i-1}$
defined by
\[R_i(\overrightarrow k)=\{\langle (v_1,\dots,v_d),(u_1,\dots,
u_d)\rangle\in\V\times\V; \exists ! r\in\{1,\dots,d\},\exists
l\in\Z,\exists j\in\{0,\dots,k_1\times\dots\times k_{i-1}-1\}\]
\begin{equation}\label{apagado11}\mbox{
such that }u_r=lk_1\dots k_i+j,v_r=k_1\dots k_{i-1}(lk_i-1)+j\mbox{ and }v_s=
u_s,\forall s\neq r\}.\end{equation}

Observe that in the simplest case, when $n=1$ and $k_1=k$, the set
of deleted bonds is precisely all nearest neighbours bonds in the
$r$-th direction where one of the endpoints has the $r$-th
coordinate multiple of $k$ and the other endpoint have the $r$-th
coordinate one unit below, for each $r=1,\dots,d$.

 Since $F^{\overrightarrow{k}}$
is subgraph of $G^{\overrightarrow{k}}$, we have the inequality
$\theta_v^{F^{\overrightarrow{k}}}(p)\leq\theta^{G^{\overrightarrow{k}}}(p),\forall v\in\V(F^{\overrightarrow{k}}).$

Now, we claim that the graphs $F^{\overrightarrow{k}}$ and
$S^{\overrightarrow{k}}$ are isomorphic. Consider the function
\begin{equation}\label{psi}
\psi:\mathbb{Z}\rightarrow \mathbb{Z}\times\prod_{i=1}^n
\{0,1,\dots,k_{n-i+1}-1\},
\end{equation}
where $$\psi(v)=\left(\Big\lfloor\frac{v}{k_1\dots
k_n}\Big\rfloor,\Big\lfloor\frac{v\mod k_1\dots k_n}{k_1\dots
k_{n-1}}\Big\rfloor,\dots,\Big\lfloor\frac{v\mod
k_1k_2}{k_1}\Big\rfloor,v\mod k_1\right).$$

Indeed the function
\begin{equation}\label{Psi}
\Psi:\mathbb{Z}^d\rightarrow
\left(\mathbb{Z}\times\Big(\prod_{i=1}^n
\{0,1,\dots,k_{n-i+1}-1\}\Big) \right)^d,\end{equation} where $\Psi
(v_1,\dots,v_d)= (\psi(v_1),\dots,\psi(v_d))$, is a graph
isomorphism between $F^{\overrightarrow{k}}$ and
$S^{\overrightarrow{k}}$. Then,
$\theta_v^{S^{\overrightarrow{k}}}(p)=\theta_{\Psi^{-1}(v)}^{F^{\overrightarrow{k}}}(p)$
for all $v\in\V(S^{\overrightarrow{k}})$, proving this lemma.

{\bf Remark.} Since
$\theta_v^{S^{\overrightarrow{k}}}(p)\leq\theta^{G^{\overrightarrow{k}}}(p),\forall v\in\V(S^{\overrightarrow{k}}),
\forall p\in [0,1]$, it holds that $p_c(G^{\overrightarrow{k}})\leq
p_c(S^{\overrightarrow{k}})$ (the second inequality in Equation
(\ref{sanduba})). Indeed, we have that the strict inequality is also
true, observing that there exists a periodic class of edges of
$G^{\overrightarrow{k}}$, which do not belong to
$S^{\overrightarrow{k}}$. (See example B in Section 3.3 of
\cite{G}).

{\it Proof of Lemma \ref{superior}}


In this proof we will use Theorem 1 of \cite{BS}, which is based on
the original idea of Campanino and Russo in \cite{CR} to prove that
the percolation threshold of the cubic lattice is bounded above by
the percolation threshold of the triangular lattice.

In the proof of Lemma \ref{inferior}, we were able to show that the
graphs $F^{\overrightarrow k}$ and $S^{\overrightarrow k}$ are
isomorphic according the function $\Psi$ defined in \ref{Psi}, where
$F^{\overrightarrow k}$ is obtained by deleting some specific edges
of $G^{\overrightarrow k}$. If we insert the respective edges again
in $S^{\overrightarrow k}$, we obtain a new graph, denoted by
$\widetilde{S}^{\overrightarrow k}$, and this latter graph is
isomorphic to $G^{\overrightarrow k}$. Formally, we have
\begin{equation}\label{sara}\widetilde{S}^{\overrightarrow k}=(\V (S^{\overrightarrow k}),\E(S^{\overrightarrow k})\cup(\displaystyle\cup_{i=1}^{n}\widetilde{R}_i(\overrightarrow{k}))).\end{equation}
Here $\widetilde{e}=\langle \widetilde{u},\widetilde{v}
\rangle\in\widetilde{R}_i(\overrightarrow {k})$ if, and only if
$\widetilde{u}=\Psi(u)$, $\widetilde{v}=\Psi(v)$ and $e=\langle u,v
\rangle \in R_i(\overrightarrow{k})$, where $R_i(\overrightarrow k)$
is defined in Equation (\ref{apagado11}).

Thus, it is enough to prove that
$$\theta^{\widetilde{S}^{\overrightarrow
k}}(p)\leq\theta^{\mathbb{Z}^{d(n+1)}}(p),\ \forall p\in[0,1].$$

For this purpose, we will show that $\widetilde{S}^{\overrightarrow
k}$ is a quotient graph of $\mathbb{Z}^{d(n+1)}$ by an automorphism
group and apply Theorem 1 of \cite{BS}.

We can write each vertex $v\in \mathbb{Z}^{d(n+1)}$ as $v=(v_1,\dots,v_d)$, with $v_i\in \mathbb{Z}^{n+1}, \forall i=1,\dots,d$. For each  $$v_j=(v_{j,1},\dots,v_{j,n+1}),$$ we will define the surjective function
\begin{equation}\label{assoc3}
\gamma:\mathbb{Z}^{n+1}\rightarrow \mathbb{Z}\times\prod_{i=1}^n\{0,1,\dots,k_{n-i+1}-1\},
\end{equation} in a recursive manner. To simplify the notation, define $(v_{j,1},\dots,v_{j,n+1})=(y_1,\dots,y_{n+1})$.

We will define $\gamma(y_1,\dots,y_{n+1})=(z_1,\dots,z_{n+1})\in\mathbb{Z}\times\prod_{i=1}^n\{0,1,\dots,k_{n-i+1}-1\}$ (here $z_1\in\Z, z_2\in\{0,\dots,k_n-1\},\dots, z_{n+1}\in\{0,\dots,k_1-1\}$), where the sequence $(z_k)_{k=1}^{n+1}$  is obtained recursively in the following way:

First, define $z_{n+1}=y_{n+1}\mod k_1$ and $t_{n+1}=\lfloor \frac{y_{n+1}}{k_1} \rfloor$.

Given $t_{i+1}$ and $z_{i+1}$ for $i=2,\dots,n$, we define $z_i$ and $t_i$ as
$$z_i=(y_i+t_{i+1})\mod k_{n+2-i} \mbox{ and } t_i=\Big\lfloor \frac{y_i+t_{i+1}}{k_{n+2-i}} \Big\rfloor .$$
Finally, define $z_1=y_1+t_2$.

Now, we define the surjection
\begin{equation}\label{assoc4}
\Gamma:\mathbb{Z}^{d(n+1)}\rightarrow (\mathbb{Z}\times\prod_{i=1}^n\{0,1,\dots,k_{n-i+1}-1\})^d,
\end{equation} where
\begin{equation}\label{sara11}
\Gamma(v)=(\gamma(v_1),\gamma(v_2),\dots,\gamma(v_d)).
\end{equation}
In words, considering the simplest case $d=1$ and $n=1$, the
function $\Gamma$ wraps $\Z^2$ onto the strip
$\Z\times\{1,\dots,k-1\}$ shifting one unit in the first coordinate
in each wind around $\Z\times\{1,\dots,k-1\}$.

Given $(z_1,\dots,z_n,z_{n+1})\in\Z^{n+1}$, for each
$j=1,\dots,n$ define the functions
$\delta_j:\Z^{n+1}\rightarrow\Z^{n+1}$ where

\[\delta_1(z_1,\dots,z_n,z_{n+1})= (z_1,\dots,z_{n-1},z_n-1,z_{n+1}+k_1),\]
\[\delta_2(z_1,\dots,z_n,z_{n+1})= (z_1,\dots,z_{n-2},z_{n-1}-1,z_n+k_2,z_{n+1}),\dots\]
\[\delta_n(z_1,\dots,z_n,z_{n+1})= (z_1-1,z_2+k_n, z_3,\dots,z_{n+1}).\]

For each $i=1,\dots,d$ and $j=1,\dots,n$ define the group
automorphism function
$\Delta_{i,j}:\Z^{d(n+1)}\rightarrow\Z^{d(n+1)}$ where
$\Delta_{i,j}(v_1,\dots,v_d)=(v_1,\dots,v_{i-1},\delta_j(v_i),v_{i+1},\dots,
v_d)$.

Using Theorem 1 of
\cite{BS} and observing that $\widetilde{S}^{\overrightarrow k}$ is
the quotient graph  $\mathbb{Z}^{d(n+1)}/\Delta$, where
$\Delta=\langle \Delta_{i,j};i=1,\dots,d,\ j=1,\dots,n \rangle$ is
the automorphism group of $\mathbb{Z}^{d(n+1)}$ generated by the set
of automorphisms $\{ \Delta_{i,j};i=1,\dots,d,\ j=1,\dots,n \}$, we have that
$\theta^{\widetilde{S}^{\overrightarrow
k}}(p)\leq\theta^{\mathbb{Z}^{d(n+1)}}(p),\ \forall p\in[0,1]$,
concluding the proof of Lemma 2.
Observe that the function $\Gamma$ defined in (\ref{sara11}) is
precisely the quotient map (the function $f$ in the proof of Theorem 1 of \cite{BS}) between $\Z^{d(n+1)}$ and
$\widetilde{S}^{\overrightarrow k}$.

\section{Final Remarks}

{\bf A)} All these results remain the same if we consider bond percolation (where each bond is open with probability $p$) instead of site percolation, since Theorem A of \cite{GM} and Theorem 1 of \cite{BS} can be used also for bond percolation.

{\bf B)} The main result, Theorem 1, can be generalized without difficulty if we consider different sequences  $\overrightarrow{k_i}=(k_{i,1},\dots, k_{i,n_i})$ for each direction $i$
(here $n_i=0$ means that there is only nearest neighbour bond in the $i$-th direction). In this case Equation \ref{eqprin} is equivalent to $$\lim_{k_{i,j}\rightarrow\infty,\forall
i,j}p_c(G^{\vec{k}})=p_c(\Z^{\sum
_{i=1}^d (1+n_i)}).$$
{\bf C)} Computational simulations (see \cite{ALS}) shows that when
$n=1$ and $d=2$ the function $\theta^k(p)$ is non decreasing in $k$,
then $p_c(G^{k+1})\leq p_c(G^k)$ as well $\lim_{k\rightarrow\infty}
p_c(G^k)=p_c(\Z^4)$ confirming Theorem 1. Based on these simulations
and in the shape of the graphs  $\widetilde{S}^{\overrightarrow k}$
we think that the following conjecture is true
\begin{conjectura}\label{conjec}
For any $p\in [0,1]$ and for any $\vec{k}$, we have $ \theta^{G^{\overrightarrow{k}}}(p)\leq\theta^{G^{\overrightarrow{k^\prime}}}(p)$, where $\overrightarrow{k^\prime}=(k_1+1,\dots,k_n+1)$.
\end{conjectura}

Indeed, it is possible to see that $G^{\overrightarrow{k}}$ is isomorphic to a quotient graph of $G^{\overrightarrow{k^\prime}}$, nevertheless $G^{\overrightarrow{k}}$ is not a quotient graph by an automorphism group, thus Theorem 1 of \cite{BS} cannot be used.

{\bf Acknowledgments.} B.N.B. de Lima and R. Sanchis are partially
supported by CNPq and FAPEMIG (Programa Pesquisador Mineiro) and
R.W.C. Silva is partially supported by FAPEMIG. We also thank the author of the careful referee report.


\begin{thebibliography}{99}




\bibitem{ALS} Atman, A.P.F., de Lima, B.N.B. and Schnabel, M.,
{\em Is the Percolation Probability on $\Z^d$ with Long Range
Connections monotone? A simulational approach.} In Preparation

\bibitem{BS}  Benjamini, I. and Schramm, O.,  {\em Percolation
beyond $\Z^d$, many questions and a few answers}. Electronic
Communications in Probability {\bf 1}, 71-82 (1996).


\bibitem{CR}  Campanino, M. and Russo, L.,  {\em An upper bound on the critical percolation probability for the
three-dimensional cubic lattice}. The Annals of Probability {\bf
13}, no. 2, 478-491 (1985).

\bibitem{FL} Friedli S., de Lima B.N.B., {\em On the truncation systems with non-Summable Interactions}, Journal
of Statistical Physics {\bf 122} 6, 1215-1236 (2006).


\bibitem{G} Grimmett G., \emph{Percolation}, 2nd edition, Springer-Verlag, Berlin, 1999.

\bibitem{GM} Grimmett,G.R. and Marstrand J.M., {\em The Supercritical Phase of Percolation is
Well Behaved}, Proc. Roy. Soc. London Ser {\bf A 430}, 439-457
(1990).





\end{thebibliography}
\end{document}